\documentstyle[12pt]{article}
\def\Bbb#1{{\mathchoice{\mbox{\bf #1}}{\mbox{\bf #1}}%
{\mbox{$\scriptstyle \bf #1$}}{\mbox{$\scriptscriptstyle \bf #1$}}}}
\def\N{\Bbb N}
\def\R{\Bbb R}
\def\C{\Bbb C}

\def\Z{\Bbb Z}
\def\T{\Bbb T}
\def\Q{\Bbb Q}
\def\E{\Bbb E}

\def\sgn{{\rm sgn}}

\def\cF{{\cal F}}

\def\sgn{{\rm sgn}}

\begin{document}
\title{Hardy martingales and Jensen's Inequality}
\author{Nakhl\'e H.\ Asmar
 and 
 Stephen J.\ Montgomery--Smith\\
Department of Mathematics\\
University of
Missouri--Columbia\\
Columbia, Missouri 65211  U.\ S.\ A.
}
\date{}
\maketitle
			     %
\begin{abstract}
\baselineskip=10 pt
Hardy martingales were introduced by Garling
and used to study analytic functions on 
the $N$-dimensional torus $\T^N$,
where analyticity is defined using a
lexicographic order on the dual group $\Z^N$.
We show how, by using basic properties of orders on
$\Z^N$, we can apply Garling's method in
the study of analytic functions on
an arbitrary compact abelian group with an arbitrary
order on its dual group.  We illustrate
our approach by giving a new and simple proof
of a famous generalized Jensen's Inequality due to
Helson and Lowdenslager \cite{hl1}.
\end{abstract}

\section{Introduction}
\newtheorem{equiv-jensen}{Lemma}[section]
Suppose that $G$ is a nonzero compact connected abelian group
with an infinite (torsion-free) dual group $\Gamma$, and
normalized Haar measure $\lambda$.
For $1\leq p<\infty$, the Banach
space of measurable functions $f$ such that
$|f|^p$ is integrable will be denoted by $L^p(G)$,
and the Banach space of essentially bounded measurable functions
on $G$ will be denoted by $L^\infty(G)$.
We 
use the symbols $\N, \Z, \R$, and $\C$ 
to denote the natural numbers, the integers, 
 the real numbers, and the complex numbers respectively.  
The circle group will be denoted by $\T$ and will be 
parametrized as $\{ e^{i t}:\ 0\leq t<2 \pi\}$.

A subset $P$ of $\Gamma$ is called an {\em order}
 if
it satisfies the following three axioms:
$$P\cap (-P)=\{0\}, P\cup (-P)=\Gamma,\ {\rm and}\ P+P=P.$$
Given an order $P\subset \Gamma$, 
we define a signum function
with respect to $P$,\ $\sgn_P$, by:
$\sgn_P(\chi)=-1,0,$\ or $1$, according as
$\chi\in (-P)\setminus\{0\}, \chi=0$,\ or
$\chi\in P\setminus \{0\}$.  
A function $f\in L^1(G)$
is called {\em analytic}
 if its Fourier transform $\widehat{f}$
vanishes off $P$.  This notion of analyticity 
was introduced by Helson and Lowdenslager \cite{hl1} and
\cite{hl2} in connection with prediction theory,
and since then it has been extensively studied
because of its independent interest.  
For $1\leq p\leq \infty$, following
\cite{hl2}, we let $H^p(G)$ 
denote
the space of analytic functions in $L^p(G)$.
To be specific about the order, we will also use the 
notation $H^p_P(G)$.

Recently,
Garling \cite{gar2} introduced Hardy martingales
and used them in \cite{gar} to prove various properties of 
analytic functions on $\T ^N$, the $N$-dimensional torus.
(See \S 3 below for a review of these notions.)
Analyticity in
\cite{gar} was defined with respect to the following 
 lexicographic order on the
dual group $\Z^N$:
\newpage
\begin{eqnarray*}
P^*	&=&	\{0\}
	\bigcup \{(m_1,m_2,\ldots,m_N)\in \Z^N:\ m_1>0,
	m_2=\ldots=m_N=0\}\\
	& &	\bigcup \{(m_1,m_2,\ldots,m_N)\in
	 \Z^N: \ m_2>0,m_3=\ldots=m_N=0\}\\
	 &&
	   \bigcup 
	 \ldots \bigcup 
		\{(m_1,m_2,\ldots,m_N)\in \Z^N:\ m_N>0\}.
\end{eqnarray*} 

Our goal in this paper is to show how Garling's approach in 
\cite{gar} can be applied in the setting of an arbitrary order
on the dual group, by using basic properties of 
orders on $\Z^N$.  We will derive these properties in Section 2.
In Section 3, we illustrate
our approach by giving a simple proof of a generalized 
Jensen's Inequality for functions
in $H^1(G)$, due to Helson and Lowdenslager \cite{hl1}.
The inequality states that for any $f\in H^1_P(G)$ we have
\begin{equation}
\left|\int_Gf(x)d\lambda(x)\right|\leq \exp \left( 
\int_G\log |f(x)|d\lambda (x)\right).
\label{jensen}
\end{equation}
This inequality plays an important role in the 
study of analytic measures
on groups, factorization of $H^1(G)$ functions, and 
the invariant subspace theory in $H^2(G)$ (see \cite{hl1} and 
\cite{hl2}).  
The only available proof
of this inequality is the original one
which uses the full strength of the methods in \cite{hl1}.  
The proof that we offer is straightforward
and follows directly from the corresponding inequality for functions
in $H^1(\T)$.

\section{Orders on $\Z^N$}
\newtheorem{remarks}{Definition}[section]
\newtheorem{defn}[remarks]{Definition}
\newtheorem{rem1}[remarks]{Remark}
\newtheorem{def2}[remarks]{Definition}
\newtheorem{rem2}[remarks]{Remarks}
\newtheorem{lemma1}[remarks]{Theorem}
\newtheorem{pure}[remarks]{Proposition}
\newtheorem{lemma2}[remarks]{Lemma}
\newtheorem{theorem3}[remarks]{Theorem}
\newtheorem{theorem4}[remarks]{Theorem}
Our goal in this section
is to prove a useful property of orders
which states that, given an arbitrary order 
$P$ on 
$\Z^N$ and a finite subset 
$E$ of $P$, there is an isomorphism
of $\Z^N$ onto itself mapping $E$ into
the lexicographic order $P^*$.
  We start by recalling from \cite[Appendix A]{hr1}
definitions and properties of 
some special subsets of discrete groups.
Throughout this section, $\Gamma$
denotes an infinite torsion-free abelian group.

\begin{remarks}
{\rm  A finite subset 
$S=\{b_1,b_2,\ldots ,b_k\}\subset \Gamma$ 
is {\em independent} if $0\not\in S$ and if 
$n_1 b_1+n_2 b_2+\ldots +n_k b_k=0$ 
implies that $n_1=n_2=\ldots =n_k=0$
(each $n_j$ is an integer).  
An infinite subset is called {\em independent} 
if every one of its finite subsets is independent.  
If an independent set $S$ generates $\Gamma$
then $S$ is called a {\em basis}.  As a convention,
if $B=\{b_1,b_2, \ldots\}$ is a basis in
$\Gamma$ and $x\in \Gamma$, we will
write
$b_j(x)$ for the $j$th coordinate of $x$ in that basis.
Hence $x=\sum_{j=1}^kb_j(x)b_j$, for some positive integer $k$.}
\label{remarks}
\end{remarks}
\begin{defn}
{\rm 
A subgroup $H$ of a torsion-free 
 group $\Gamma$ is called {\em pure} in $\Gamma$ 
 if whenever $x\in \Gamma , n\neq 0,\ n x \in H$, then
 $x\in H$.  }
\label{defn}
\end{defn}
Pure subgroups will arise in our
proofs as kernels of 
homomorphisms.  The concept of pure subgroups
is very important in the study of the structure of abelian
groups.  (See \cite[Appendix A]{hr1}.)
In what follows, we list some basic properties
 related to this concept, which will be needed in the sequel.
\begin{rem1}
{\rm 
 (a) Suppose that $H$ and $K$
are subgroups of $\Gamma$, and that
$H\subset K\subset \Gamma$.  It is easy to see that if $K$ 
 is pure in $\Gamma$ and $H$ is pure 
 in $K$, then $H$ is pure in $\Gamma$.  \\
The following is a very important property of 
pure subgroups of $\Z^N$.\\
(b) If $H$ is a proper pure subgroup of $\Z^N$,  
then $H$ is isomorphic to $\Z^\nu$ for some $\nu <N$.  
In this case, there is a basis of $\Z^N$ of the form 
 $\{b_1,b_2,\ldots,b_\nu,\ldots,b_N\}$ such that 
  $\{b_1,b_2,\ldots,b_\nu\}$  is a basis of $H$ 
(see \cite[Theorem A. 26]{hr1}).
The number $\nu$ is called the {\em rank},
or the {\rm dimension}, of $H$, and will be denoted
$r(H)$.  
\\
Using (a), we can generalize (b) as follows.\\
(c)  If 
  $\{0\}\subset H_1
  \subset H_2\subset\ldots \subset H_k=\Z^N$ 
  is a sequence of subgroups of $\Z^N$ 
  such that $H_j$ is a proper pure subgroup of 
  $H_{j+1}$, then there is a basis of $\Z^N$, $B=\{b_1,\ldots,b_{r(H_1)},\ldots,b_{r(H_j)},\dots,b_{r(H_k)}\}$,\ 
  such that, for $j=1,\ldots,k$,
  $\{b_1,\ldots,b_{r(H_j)}\}$ 
  is a basis of $H_j$.\\
}
\label{rem1}
\end{rem1}
In addition to these notions, we need 
the following result, \cite[Theorem (2.5)]{asm1}, that describes
orders in $\Z^N$ in terms of a 
decreasing sequence of subgroups
and corresponding separating real-valued homomorphisms.

\begin{lemma1}
Let $P$ be an arbitrary order on $\Z^N$.  There are a 
strictly increasing sequence of subgroups
\begin{equation}
\{0\}=C_0 \subset C_1\subset 
\ldots \subset C_{k-1}\subset C_k=\Z^N,
\label{finite-chain}
\end{equation}
and a sequence 
$\{L_j\}_{j=1}^k$ 
of real-valued homomorphisms of $\Z^N$
such that, for $j=1,2,\ldots,k$, we have \\
(i)  $ L_j\left(C_{j-1}\right)=\{0\}$;\\
(ii)  $\sgn_P(\chi)=\sgn (L_j(\chi))$,
for all $\chi \in C_j\setminus C_{j-1}$.  
\label{lemma1}
\end{lemma1}
For the lexicographic order $P^*$
(see Section 1), Theorem \ref{lemma1} is obvious.  
In this case, 
we have $k=N$; $C_0=\{0\}, 
C_1=\{x\in \Z^N:\ x_2=\ldots=x_N=0\},\
\ldots,\  
C_{N-1}=\{x\in \Z^N:\ x_N=0\},\ 
C_N= \Z^N;$ and 
$L_j(x_1,x_2,\ldots,x_N)=x_j$ for $j=1,\ldots,N$.\\

Since $C_{k-1}$ is the kernel of a
homomorphism of $\Z^N$, it follows immediately that
$C_{k-1}$ is a pure subgroup of $\Z^N$.  Similarly,
$C_{k-2}$ is a pure
subgroup of $C_{k-1}$.  Now by Remark \ref{rem1} (a),
it follows that $C_{k-2}$ is a pure subgroup of $\Z^N$.
Continuing in this fashion, we obtain the following simple
proposition. 
\begin{pure}
In the notation of Theorem \ref{lemma1},
we have that every $C_j$, $j=1,\ldots,k$, is a
pure subgroup of $\Z^N$.
\label{pure}
\end{pure}
Let $H$ 
denote a nonzero pure subgroup of $\Z^N$.  
Write an element $x\in H$
as $x=\sum_{l=1}^{r(H)} b_l(x) b_l$ where 
$\{b_l\}_{l=1}^{r(H)}$ is a basis of $H$.  
A nonzero 
homomorphism $L$ of $H$ into $\R$ is said to {\em have 
integer coefficients} if 
$$L(x)=  \sum_{l=1}^{r(H)} \alpha_l b_l(x) $$
where $\{\alpha_l\}_{l=1}^{r(H)}\subset \Z$.  
Note that if $L$ has integer coefficients with respect to one basis,
then it has integer coefficients with respect to all bases.
\begin{lemma2}
Suppose that $H\neq \{0\}$ is a pure
subgroup of $\Z^N$, and $L$ is a
nonzero homomorphism of $H$ 
with integer coefficients.  
Let $K=\ker L=\{x\in H:\ L(x)=0\}$.  Then 
there is an element $h\neq 0 \in H$ such that 
$$H=K\oplus <h>,$$
where $<h>$ denotes the subgroup of $\Z^N$
generated by $h$.
In other words, $\ker L$ has codimension 1 in $H$.
\label{lemma2}
\end{lemma2}
{\bf Proof.}  Since $H$ is a pure
subgroup, we may without loss of generality assume that
$H\equiv \Z^\nu$.  The homomorphism $L$
may be considered as a homomorphism from the linear 
vector space $\Q^\nu$ (over the field $\Q$)
into $\Q$.  Since the rank of this mapping is
clearly 1, its kernel has dimension $\nu -1$,
and so contains $\nu -1$ independent vectors over
$\Q$.  By multiplying these vectors by
a large enough integer, we get $\nu -1$ independent vectors
of $\Z^\nu$ belonging to $\ker L$.  These vectors
necessarily form a basis for $K$.
Now, since 
$\ker L$ is a pure subgroup of 
$\Z^\nu$, the lemma follows from Remark \ref{rem1} (c).
\\

Some more items of notation are
needed before we state the main theorem of this section.
To be specific about an order $P$ on a group $\Gamma$,
we will sometimes write $(\Gamma, P)$.   
Given a nonvoid finite subset $S$ of $(\Z^N,P)$, 
we will write
$S_j$ to denote the set
$$S_j=S\cap (C_j\setminus C_{j-1}),\ \ j=1,2,\ldots,\ k,$$
where $C_j$ is the subgroup of $\Z^N$ given by Theorem 
\ref{lemma1}.  According to Theorem 
\ref{lemma1}, there are at most $k\leq N$ such sets $S_j$,
uniquely determined by $P$.
The case of a lexicographic order is particularly interesting to us.
In that case, we will use the notation $S^*_j$ for the sets
$S_j$.  Hence, for $j=1,\ldots, N$, we have
$$S_j^*=\{x\in S:\ x_j\neq 0, x_{j+1}=\ldots=x_N=0\}.$$
\begin{theorem4}
Let $S$ be a finite nonvoid subset of $\Z^N$,
let $P$ be an arbitrary order on $\Z^N$,
and let $\{C_j\}_{j=0}^k$ be as in Theorem \ref{lemma1}.  
There is an isomorphism
$$\psi\ :(\Z^N,P)\longrightarrow  (\Z^N,P^*)$$
such that $\psi(S\cap P)\subset P^*$, and
$\psi(S\cap (-P))\subset (-P^*)$.  Moreover,
for $j=1,\ldots,k$, we have
$$\psi\left(S_j\right)=\left(\psi(S)\right)^*_{r(C_j)}.$$ 
\label{theorem4}
\end{theorem4}
{\bf Proof.} We apply Theorem \ref{lemma1} and use its notation.
We will construct a basis $B$ of $\Z^N$ of the form
$$B=\{c_1,\ldots,c_{r(C_1)},\ldots,c_{r(C_j)},\ldots,c_{r(C_k)}\},$$
so that, in that basis, every
element of $S_j$ has a nonzero
$r(C_j)$-th component, and $x\in S_j\cap P$ if and
only if $x\in S_j$ and $c_{r(C_j)}(x)>0$.  Then
the theorem will follow by setting
$$\psi\left(
\sum_{l=1}^N c_l(x)c_l\right)=\sum_{l=1}^N c_l(x)e_l
$$
where $\{e_1,\ldots,e_N\}$ denotes the standard 
basis in $\Z^N$.  
We now proceed to show how to construct $B$.
Without loss of generality,
we may assume that $S_j$ is not empty for all $j=1,\ldots,k$.
Let
$$B_1=\{b_1,\ldots,b_{r(C_1)},\ldots,b_{r(C_j)},\ldots,b_{r(C_k)}\},$$
denote a basis of $\Z^N$ with the property that
$\{b_1,\ldots,b_{r(C_1)},\ldots,b_{r(C_j)}\}$
is a basis for $C_j$ (Remark \ref{rem1} (c)).
For $x\in \Z^N$, we write  
$x=\sum_{l=1}^Nb_l(x)b_l$.
Expressing the homomorphism
$L_j$ in the basis $B_1$, we have for
 $x\in C_j$, 
\begin{equation}
L_j(x)=\sum_{l=r(C_{j-1})+1}^{r(C_j)}
\beta_lb_l(x),
\label{lj}
\end{equation}
where $\beta_l\in\R$, because $L_j (C_{j-1})=\{0\}$.  Since,
by Theorem \ref{lemma1} (ii),
$L_j(x)>0$ for all $x\in S_j\cap P$,
and $L_j(x)<0$ for all $x\in S_j\cap (-P)$, 
and since $S$ is finite,
we can replace the coefficients $\beta_l$ in (\ref{lj})
by integers $\alpha_l$ so that
$$\sum_{l=r(C_{j-1})+1}^{r(C_j)}
\alpha_l b_l(x)>0$$
if $x\in S_j\cap P$, and 
$$\sum_{l=r(C_{j-1})+1}^{r(C_j)}
\alpha_l b_l(x)<0$$
if $x\in S_j\cap (-P)$.
(First replace the real numbers
$\alpha_l$ in (\ref{lj}) by rational numbers,
then multiply by a sufficiently large positive integer.) 
Define a homomorphism $L_j^*$ for all $x\in \Z^N$ by
$$L_j^*(x)=\sum_{l=r(C_{j-1})+1}^{r(C_j)}
\alpha_l b_l(x).$$
Plainly, $L^*_j(x)\neq 0$, for all $x\in S_j$.
Let 
$D_j=\{x\in C_j:\ L_j^*(x)=0\}$.
Then $D_j\supset C_{j-1}$, and,
by Lemma \ref{lemma2}, $D_j$ has codimension 1
in $C_j$.  Let $h_j\in C_j$ be such that
$$C_j=D_j\oplus <h_j>.$$
We can and do choose $h_j$ so that
$L_j^*(h_j)>0$.  Now consider the basis
$$B=\{c_1,\ldots,c_{r(D_1)},h_1=c_{r(C_1)},
\ldots,  
c_{r(D_j)},h_j=c_{r(C_j)},
\ldots,  
c_{r(D_k)},h_k=c_{r(C_k)}
\},$$
where
$\{
c_1,\ldots,c_{r(D_1)},h_1,
\ldots,  
c_{r(D_j)}
\}$
is a basis for $D_j$, and hence 
$$\{
c_1,\ldots,c_{r(D_1)},h_1,
\ldots,  
c_{r(D_j)}, h_j
\}$$
is a basis for $C_j$.
For $x\in S_j$ we have
$L^*_j(x)\neq 0$, and so $h_j(x)\neq 0$.  Also,
  $x\in S_j\cap P$
if and only if $x\in S_j$ and $L^*_j(x)>0$,
if and only if 
$x\in S_j$ and $h_j(x)>0$.
This shows that $B$ has the desired properties and
completes the proof of the theorem.
\section{Hardy martingales and Jensen's Inequality}
\newtheorem{prop1}{Proposition}[section]
\newtheorem{prop2}[prop1]{Proposition}
\newtheorem{density}[prop1]{Lemma}
We start this section by reviewing the concept of Hardy martingales
from \cite{gar}.
Let $N$ be a fixed positive integer, and let
$e^{i\theta_n}$ denote the $n$-th coordinate evaluation function on $\T^N$.
Let $\cF_n=\sigma(e^{i\theta_1},e^{i\theta_2},\ldots,e^{i\theta_n})$
denote the $\sigma$-algebra generated by the first $n$ coordinate functions.
For $f\in L^1(\T^N)$, the conditional expectation of $f$
with respect to $\cF_n$ will be denoted 
$\E(f|\cF_n)$.  For $n=0,1,2,\ldots,N$, the function 
$\E(f|\cF_n)$
is constructed from $f$  
by projecting the Fourier transform of $f$ 
on $\Z^n$, where here $\Z^0=\{0\}$, and
$$\Z^n=\{(k_1,k_2,\ldots,k_n,0,\ldots,0):\ k_j\in \Z\}.$$
The finite sequence 
$\left(  \E(f|\cF_n) \right)_{n=1}^N$ forms a martingale relative to  $\left(\cF_n\right)_{n=1}^N$.
(A detailed analysis of martingales on groups
with a decreasing sequence of subgroups
can be found in \cite[Chapter 5]{eg}.)
A complex-valued martingale $(f_n)$
on $\T^N$ is called a {\em Hardy martingale} if 
$\E(f_{n+1}e^{ik \theta_{n+1}}|\cF_n)=0$ for $k>0$ and all $n=0,\ldots,N-1$.

Hardy martingales arise naturally when studying
analytic functions in $H^1(\T^N)$.
Indeed, as observed in \cite{gar}, for $1\leq p\leq \infty$, 
we have
\begin{equation}
H^p_{P^*}(\T^N)=\left\{f\in L^p(\T^N):\ (\E(f|\cF_n))_{n=1}^N
\ {\rm is\ a \ Hardy\ martingale } \right\},
\label{hp*}
\end{equation}
where here, as in \S 1, 
$P^*$ denotes the lexicographic order on $\Z^N$.

It is instructive to justify this fact,
and show, in the process,
how Hardy martingales are related to the usual 
Hardy spaces on the circle group.
For this purpose, we recall the notion 
of martingale differences series.
Let
\begin{equation}
d_0(f)=\int_{\T^N} fdx,
\label{d0}
\end{equation}
and for $j=1,\ldots, N$, let
$$d_j(f)=\E(f|\cF_j)-\E(f|\cF_{j-1}).$$
Thus $d_j (f)$ ($j=1,\ldots,N$)
is constructed from $f$ by projecting the Fourier transform
of $f$ on the set difference $\Z^j\setminus \Z^{j-1}$.
So
$d_j(f)$\ may be formally represented
as the Fourier series:
\begin{equation}
d_j(f)=\sum_{k=-\infty,k\neq 0}^\infty 
f_{j,k}(\theta_1,\ldots,\theta_{j-1})
e^{i k \theta_j}, 
\label{dj1}
\end{equation}
 where $ f_{j,k}(\theta_1,\ldots,\theta_{j-1})$ is a function 
 of $\theta_1,\ldots,\theta_{j-1}$ only. 
 We have the {\em martingale difference series}
decomposition
\begin{equation}
f=\sum_{j=0}^N d_j(f).
\label{martingale-difference}
\end{equation}
As observed in \cite{gar}, and is easy to check,
$f$ is in $H^1_{P^*}(\T^N)$ if and only if,
for $j=1,\ldots,N$,
\begin{equation}
d_j(f)=\sum_{k=1}^\infty f_{j,k}(\theta_1,\ldots,\theta_{j-1})
e^{i k \theta_j}. 
\label{dj}
\end{equation}
 Hence, as a function of $\theta_j$, the 
 function $d_j(f)$ belongs to $H^1(\T)$ for $j=1,2,\ldots,N$.
 From these observations, (\ref{hp*}) 
follows easily.  

We now return to  
inequality (\ref{jensen}), and prove a special
case of it.
\begin{prop1}
Suppose that $f\in H^1_{P^*}(\T^N)$.  Then,
$$ \exp \left(\int_{\T^N}\log | \sum_{j=0}^n d_j(f) |dx\right)
   \le
   \exp \left(\int_{\T^N}\log | \sum_{j=0}^{n+1} d_j(f) |dx\right)$$
and
$$\left|\int_{\T^N}fdx\right|\leq 
	\exp \left(\int_{\T^N}\log |f|dx\right).$$
\label{prop1}
\end{prop1}
{\bf Proof.}  
The second inequality follows by applying the first inequality 
repeatedly (for $n=0,1,\ldots,N-1$) and using 
(\ref{d0}) and (\ref{martingale-difference}).
Let us prove the first inequality.  
As we observed above, the function
 $\theta_{n+1}\mapsto d_{n+1}(f)$ is in $H^1(\T)$.
Using the one-dimensional version of Jensen's Inequality
(\ref{jensen}) for functions in $H^1(\T)$ 
(see \cite[Inequality (3.2), and Theorem 3.11]{ka}), 
we obtain
$$ 
\log \left|\frac1{2\pi} \int_{\theta_{n+1}}
\left(\sum_{j=0}^{n+1}d_j(f)
\right) d\theta_{n+1}\right|
\leq
\frac1{2\pi} \int_{\theta_{n+1}}
\log \left|\sum_{j=0}^{n+1}d_j(f)\right| d\theta_{n+1}.$$
But 
$$
\int_{\theta_{n+1}}
d_{n+1}(f) d\theta_{n+1}=0,
$$
and $d_j(f)$ is constant in $\theta_{n+1}$ for $j=0,1,\ldots,n$.  
Thus,
$$ 
\log \left|
\sum_{j=0}^{n}d_j(f)
 \right|
\leq
\frac1{2\pi} \int_{\theta_{n+1}}
\log \left|\sum_{j=0}^{n+1}d_j(f)\right| d\theta_{n+1}.
$$
Integrating with respect to the remaining variables, we get
$$
\int_{\T^N}\log \left| \sum_{j=0}^{n} d_j(f) \right| dx
\leq
\int_{\T^N}\log \left| \sum_{j=0}^{n+1} d_j(f) \right| dx,
$$
which completes the proof of the proposition.

For the remainder of the proof, we need the following density result.
\begin{density}
Let $G$\ be a compact abelian group with dual ordered by
$P$, and let $Y$\ be a dense subspace of $H^1_P(G)$.
If (\ref{jensen}) is true for all $f$\ in $Y$, then
(\ref{jensen}) is true for all $f$\ in $H^1_P(G)$.
\label{density}
\end{density}
{\bf Proof.}  
We first note that for a given function $f \in L^1(G)$\ that (\ref{jensen})
holds if and only if for all $0<p<1$\ we have
\begin{equation}
   \left| \int_G f d\lambda\right|^p \le
   \int_G |f|^p d\lambda .
\label{pjensen}
\end{equation}
Indeed (\ref{jensen}) follows from (\ref{pjensen}) 
by letting $p$ tend to zero
(see \cite[(13.32)(ii)]{hs}).  Now suppose that (\ref{jensen}) holds.
Then, for any $0<p<1$, we have
\begin{eqnarray*}
\left| \int_G f d\lambda\right|^p 
					&\leq&
\exp\left(\int_G \log\left(|f|^p\right) d\lambda \right)
					\leq
\int_G |f|^p d\lambda,
\end{eqnarray*}
where the last inequality follows from 
\cite[(13.32) (i)]{hs}.

Now fix $0<p<1$\ and $f \in H^1_P(G)$.  
Let $(f_n)$\ be a sequence in $Y$\
such that $f_n \to f$\ in $L^1(G)$.  
We have that 
$f_n \to f$ in $L^p$, and from the inequality
$\int_G\left| |f_n|^p-|f|^p\right|d \lambda\leq 
\int_G\left|f_n-f\right|^p d\lambda$, it follows that
$\int_G\left|f_n\right|^p d\lambda\to 
\int_G\left|f\right|^p d\lambda$
(see \cite[Theorem (13.17) and (13.25) (a)]{hs}).
Since (\ref{pjensen}) holds for every $f_n$, it
follows immediately that it also holds for $f$.

The next step is to establish (\ref{jensen})
for arbitrary orders on $\Z^N$.  At this point we will
appeal to Theorem \ref{theorem4}.
\begin{prop2}
Suppose that $f\in H^1_{P}(\T^N)$.  Then, 
$$\left|\int_{\T^N}fdx\right|\leq \exp\left(\int_{\T^N}\log|f|dx\right).$$
\label{prop2}
\end{prop2}
{\bf Proof.}  
Lemma \ref{density} shows that it is enough
to consider the case of a trigonometric polynomial 
$f\in H^1_P(\T^N)$.  Write 
$$f=\sum_{\chi\in S}a_\chi \chi$$
where $a_\chi\in \C$, and $S$ is a finite subset of $P$.
Apply Theorem \ref{theorem4} and obtain an isomorphism $\psi$ 
of $\Z^N$
such that $\psi(S)\subset P^*$.  
Let $\phi$ be the adjoint homomorphism of $\psi$.
Thus $\phi$ is an automorphism
of $\T^N$ onto itself such that
$$\psi(\chi)(x)=\chi\circ\phi(x)$$
for all $\chi\in \Z^N$ and all $x\in \T^N$.  Moreover,
for all $\chi\in S$, the character
$\chi\circ\phi=\psi(\chi)$ is in $P^*$.
Since $\phi$ is an automorphism of $\T^N$, it maps 
the 
normalized Haar measure to itself.
Using this last observation and Proposition \ref{prop1}, we obtain
\begin{eqnarray*}
\log\left|\int_{\T^N}f\, dx\right|
				&=&
\log\left|\int_{\T^N}f\circ\phi\, dx\right|
				=
\log\left|\int_{\T^N}\sum_{\chi\in S}a_\chi \psi(\chi)\, dx\right|\\
				&\leq&
\int_{\T^N}\log\left|\sum_{\chi\in S}a_\chi \psi(\chi)\right| dx
				=
\int_{\T^N}\log\left|f\right| dx,				
\end{eqnarray*}
which completes the 
proof of the proposition.

We are only a step away from completing the proof of (\ref{jensen}).
The standard reduction that remains to be done
is based on the Weil formula.  We present the details for the sake of 
completeness.

{\bf Proof of (\ref{jensen}).}  
Throughout this proof, $P$ will denote an arbitrary order on 
$\Gamma$.  By Lemma \ref{density}, 
it is enough to consider 
trigonometric polynomials in $H^1_P(G)$.  Let
$$f=\sum_{\chi\in S}a_\chi \chi,$$
where $S$ is a nonvoid 
finite subset of $P$, and $a_\chi\in \C$.  
Let $\langle S \rangle$ denote the subgroup of $\Gamma$
generated by $S$.  Since $\Gamma$ is torsion-free,
$\langle S \rangle$ is isomorphic to $\Z^N$ for some positive integer $N$. 
Let $G_0$ denote the annihilator in $G$ of $\Z^N$.
Then $G/G_0$ is topologically isomorphic to $\T^N$
and its dual group is $\langle S \rangle\equiv \Z^N$.  We order
$\Z^N$ by intersecting it with $P$.
Let $\Pi$ denote the
natural homomorphism of $G$ onto 
$G/G_0$.
Since $f$ is constant on cosets of $G_0$,
there is a trigonometric polynomial on $\T^N$, $f^\dagger$,
such that
$$f=f^\dagger \circ \Pi.$$
Clearly $f^\dagger\in H^1_{\Z^N\cap P}(\T^N)$.
Now, using the Weil formula \cite[Theorem (28.54) (iii)]{hr2} 
and Proposition \ref{prop2},
we find that
\begin{eqnarray*}
\log\left|\int_G fd\lambda\right|
			&=&
\log\left|\int_{\T^N} f^\dagger dx\right| \\
			&\leq&
\int_{\T^N} \log\left|f^\dagger \right| dx
			=
\int_G \log\left|f \right| d\lambda,
\end{eqnarray*}
which yields the desired inequality and completes the proof of (\ref{jensen}).


{\bf Acknowledgements}  The work of the authors was supported
by separate grants from the National Science Foundation (U.S.A.).

\end{document}